\documentclass{scspaperproc}

\usepackage{latexsym}
\usepackage{graphicx}

\usepackage[pdftex,colorlinks=true,urlcolor=blue,citecolor=black,anchorcolor=black,linkcolor=black,bookmarks=false]{hyperref}
\usepackage[english]{babel}
\usepackage{amsmath}
\usepackage{array}
\usepackage{amssymb}
\usepackage{multirow}
\usepackage{float}
\usepackage{longtable}
\usepackage{hyperref}
\usepackage{cleveref}
\usepackage{graphicx}
\usepackage{caption}
\usepackage{subcaption}
\usepackage{hyphenat}
\IfFileExists{accessibility.sty}{\usepackage{accessibility}}{}

\makeatletter
\@ifpackageloaded{hyperref}{%
  \AtBeginDocument{\hypersetup{
    colorlinks=true,
    linkcolor=black,
    citecolor=black,
    urlcolor=blue,
    pdfborder={0 0 0}
  }}%
}{%
  \usepackage[
    colorlinks=true,
    linkcolor=black,
    citecolor=black,
    urlcolor=blue,
    pdfborder={0 0 0}
  ]{hyperref}%
}
\makeatother

\hypersetup{
  pdftitle={ANNSIM 2026 Proceedings Paper},
  pdfauthor={Proceedings Authors},  
  pdfsubject={ANNSIM 2026 Conference},
  pdfkeywords={}
}

\begin{document}

%
%
\SCSpagesetup{Pang, Delahaye, and Clarke}

\def\SCSconferencename{Annual Modeling and Simulation Conference}

\def\SCSconferenceacro{ANNSIM'26}

\def\SCSpublicationyear{2026}

\def\SCSconferenceeditors{G. Rabadi, V. Prabhu, R. Cárdenas, A. Bany Abdelnabi, J. Jabbour, and M. Germanos}

\def\SCSconferencedates{May 4-7}

\def\SCSconferencevenue{University of Central Florida, Orlando, Florida, USA}

\title{Geometric Trajectory Optimization for TRACON Arrivals: An NLP Approach with ATC Vectoring Maneuver Modeling}

\author[\authorrefmark{1}]{Yutian Pang}
\author[\authorrefmark{2}]{Daniel Delahaye}
\author[\authorrefmark{1}]{John-Paul Clarke}

\affil[\authorrefmark{1}]{Department of Aerospace Engineering and Engineering Mechanics, The University of Texas at Austin, USA}
\affil[\authorrefmark{2}]{OPTIM, École Nationale de l’Aviation Civile, France}
\maketitle

\section*{Abstract}
Terminal airspace congestion remains a major bottleneck in the global air traffic network. Although the Aircraft Sequencing and Scheduling Problem (ASSP) has been widely studied, many methods rely on simplified node-link abstractions that ignore the practical flight path, producing schedules that can be hard to execute under real airspace geometric constraints. This paper introduces a high-fidelity trajectory optimization framework for Terminal Arrival Sequencing and Scheduling (TASS) that explicitly models controller vectoring maneuvers. We formulate a single-stage nonlinear programming (NLP) model with a weighted objective function that optimizes Baseleg path extension and segment-wise speed profiles for arriving aircraft. The model enforces nonlinear geometric coupling between Baseleg extensions and the required Radius-to-Fix (RF) turn for Final Approach Fix (FAF) intercept through closed-form expressions. Under the First-Come-First-Served (FCFS) rule, it generates separated, operationally feasible trajectories aligned with controller practice. Monte Carlo simulations with ample runs on the A80 TRACON in Atlanta demonstrate minimal separation violations below the runway capacity threshold through tactical path stretching, with violations occurring only when arrival rates exceed the maximum tolerable rate.

\textbf{Keywords:} TRACON Automation, Terminal Operation, Arrival Sequencing, Landing Separation, Optimization

\section{Introduction}\label{sec: intro}
The sustained growth in global air traffic demand has placed immense pressure on the National Airspace System (NAS), with the Terminal Maneuvering Area (TMA) emerging as the primary bottleneck restricting system-wide throughput. In complex metroplex environments, such as the Atlanta and Los Angeles regions \cite{ren2009contrast}, the interdependence of arrival and departure flows requires precise synchronization to maintain safety standards while maximizing runway utilization. As traffic density increases, the traditional First-Come-First-Served (FCFS) operational paradigm becomes increasingly inefficient, leading to excessive fuel burn, noise emissions, and unrecoverable delays \cite{ng2024optimization}. Consequently, the Aircraft Sequencing and Scheduling Problem (ASSP) has become a focal point of research, necessitating advanced decision support tools that can transition operations from reactive control to trajectory-based strategic planning.

Early efforts to automate terminal operations, such as the Center-TRACON Automation System (CTAS) \cite{denery1995center} and the Final Approach Spacing Tool (FAST) \cite{halverson1992systems}, successfully introduced the concept of time-based metering. These systems significantly improved controller situational awareness and reduced delays by providing advisory sequencing. However, as noted in another literature \cite{diffenderfer2013automated}, existing procedures often rely on static intervals and manual coordination between the Tower and TRACON to manage gaps, which can compromise throughput when arrival and departure streams interact. Furthermore, while these systems provide estimated times of arrival, they often lack the high-fidelity trajectory optimization required to automatically generate the specific speed and path controls necessary to achieve those times precisely.

The academic literature has largely addressed the ASSP through combinatorial sequencing and trajectory control. On the sequencing side, the problem is frequently formulated as a Mixed-Integer Linear Program (MILP). For instance, \cite{desai2016optimization} demonstrated that optimizing flight sequences over the entire TMA region, rather than just the runway threshold, can reduce delays by nearly 30\% compared to FCFS. However, the ASSP is known to be NP-hard, and adding realistic operational constraints often renders exact MILP formulations computationally intractable for real-time applications. This has led researchers to favor metaheuristic approaches (i.e., genetic algorithms or simulated annealing \cite{gui2025metaheuristic, dhief2023meta}) to find optimal solutions within limited time requirements. While effective at handling binary sequencing variables, these heuristic methods often simplify aircraft flight mechanics, yielding sequences that are theoretically efficient but operationally unsafe.
\begin{figure}
    \centering
    \includegraphics[width=0.75\textwidth]{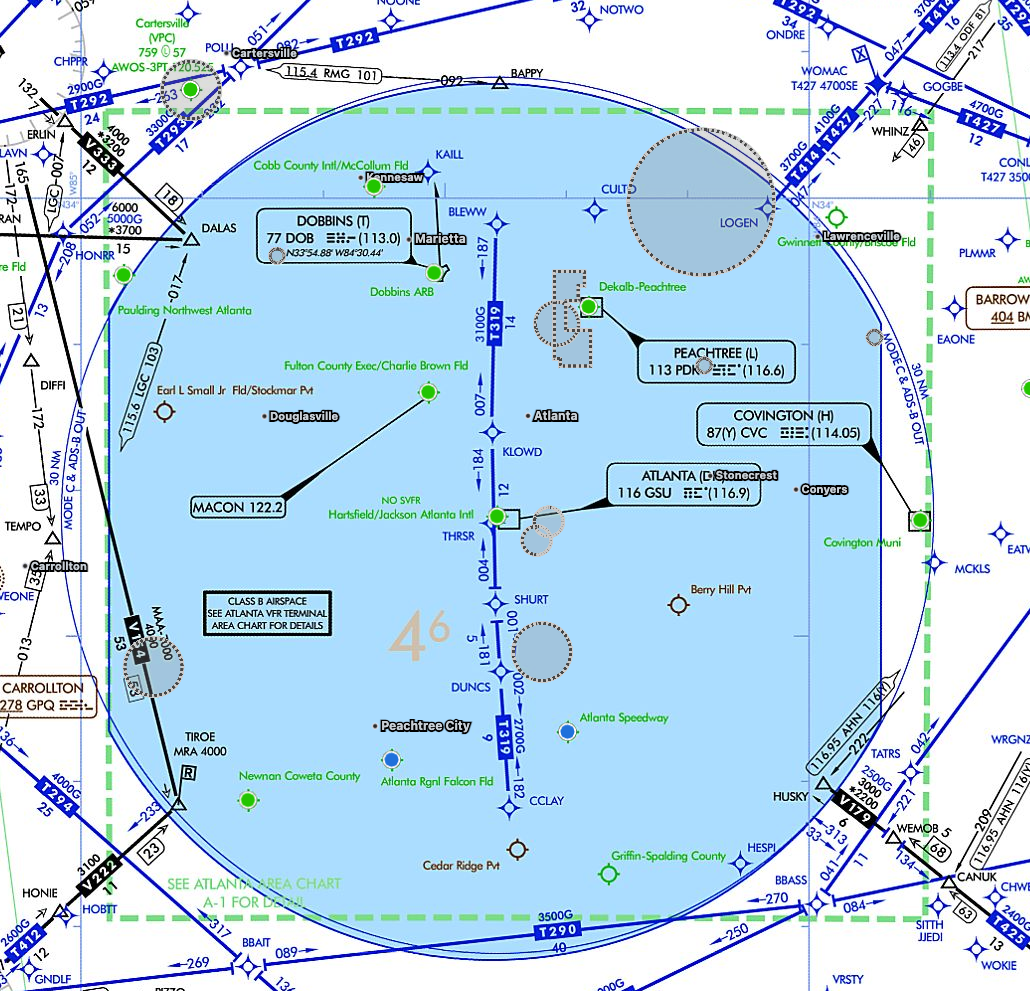}
    \caption{Entry fixes for the A80 ATL sector. Aircraft enter the terminal area via the entry fix of each feeder gates: DALAS (northwest), LOGEN (northeast), HUSKY (southeast), and TIROE (southwest). Figure courtesy of ForeFlight \cite{ForeFlight2025}.}
    \label{fig: A80_entry_fixes}
\end{figure}

Conversely, trajectory-based studies emphasize the continuous control of aircraft to absorb delay. \cite{liang2017integrated} proposed integrated sequencing using advanced avionics and point merge structures, while \cite{dhief2023meta} explored automated vectoring. However, a critical gap remains in the literature regarding the mathematical modeling of the standard \textit{trombone} pattern (Baseleg extension) used in high-density TRACONs. Most existing optimization frameworks approximate delay absorption as a purely temporal variable or assume fixed path lengths, neglecting the complex geometric coupling between the downwind extension length ($d_i$) and the subsequent Radius-to-Fix (RF) turn to the Final Approach Fix (FAF). Ignoring these geometric constraints can lead to solutions that require impossible turn radii or violate airspace sector boundaries.

To address these limitations, this paper proposes a high-fidelity trajectory optimization framework that explicitly models the operational geometry of terminal arrival flows. We adopt the A80 ATL sector as the simulated environment, utilizing the four standard feeder gates (DALAS, LOGEN, HUSKY, TIROE) \cite{ZTL_A80_LOA_2022}. The primary contribution of this work is the formulation of a Nonlinear Programming (NLP) model that mathematically formalizes the \textit{trombone} maneuver. By treating the path extension $d_i$ as a continuous decision variable, we derive closed-form expressions for the resulting RF turn geometry, ensuring that every generated solution satisfies the kinematic constraints of the aircraft. We co-optimize this geometric extension with segment-wise speed profiles under a FCFS sequence. This approach bridges the gap between abstract scheduling algorithms and operational reality, providing a benchmark for how much delay can be physically absorbed through path stretching and speed control before holding patterns become necessary.

\section{Methodology}\label{sec: method}
\begin{figure}[t]
    \centering
    \includegraphics[width=\textwidth]
    {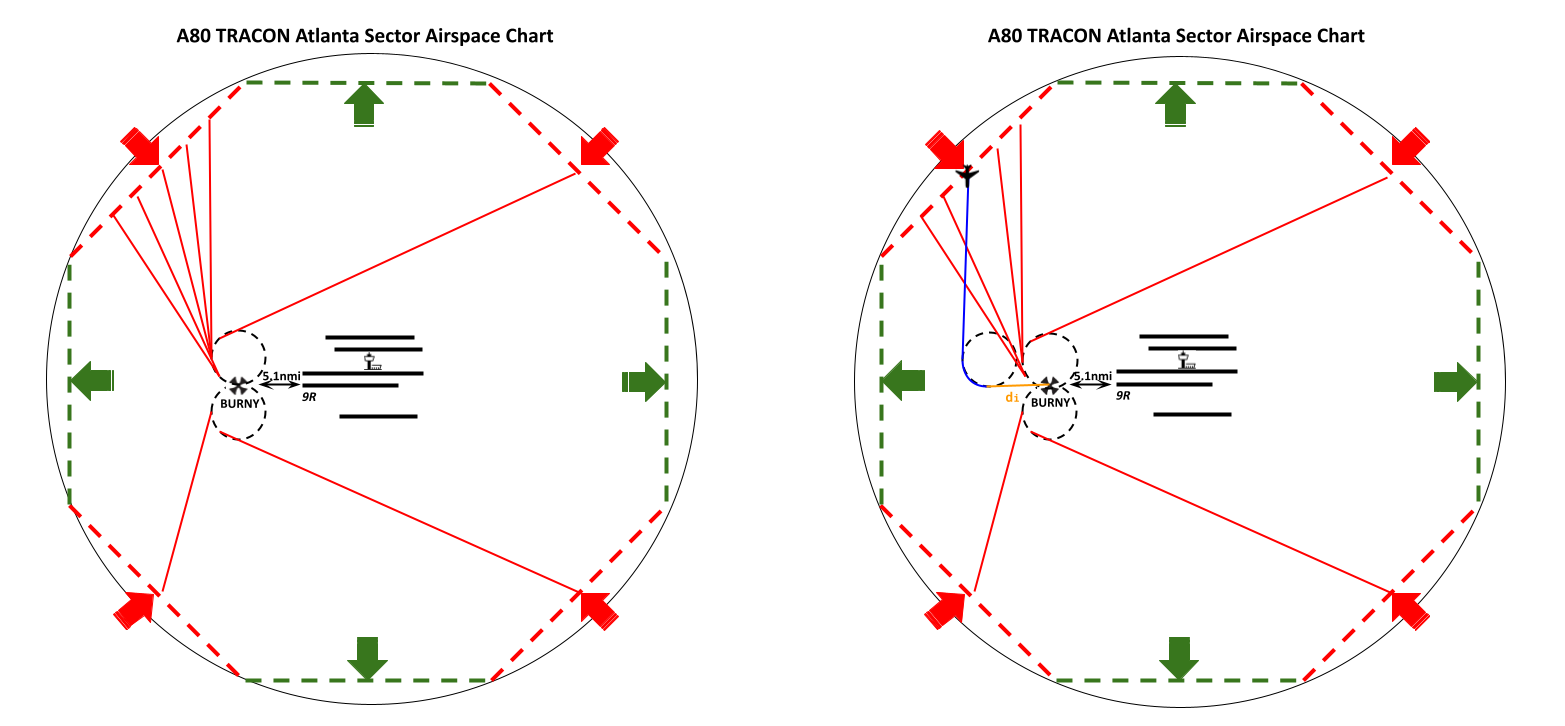}
    \caption{Schematic diagram of the Atlanta VORTAC 30 DME (TCP) of the Atlanta Large TRACON (A80), where A80 traffic controller takes over the aircraft from the Atlanta Air Route Traffic Control Center (ZTL). The circular airspace approximately represents the A80 ATL Sector approach control boundary, with arrival procedures from four feeder gates (colored in red). In this study, the arrivals are looking for a west landing of runway 9R at KATL, where the Final Approach Fix (FAF) is BURNY. In this figure, left panel shows the ideal case where aircraft arriving from each corner each takes their shortest path and an radial-to-fix (RF) turn to FAF. However, in real practice, as shown in the right panel, the arrival aircraft might have to turn to an extension of FAF (i.e., $d_i$), to ensure proper aircraft separation at FAF.}
    \label{fig: A80}
\end{figure}

In this section, we present our proposed terminal vectoring model for strategic arrival traffic management. As illustrated in \Cref{fig: A80}, the fundamental control mechanism modeled is the \textit{trombone} maneuver during ATC vectoring, where aircraft extend the Baseleg by a variable distance $d_i$ upstream of the FAF. This extension allows for precise temporal spacing without the need for fuel-inefficient holding patterns.
We formulate this problem as a continuous NLP model. The objective is to determine the optimal path extension $d_i$ and the associated segment speeds (tangent leg speed $v_{L,i}$, RF turn speed $v_{\theta,i}$, and final approach speed $v_{f,i}$) for each aircraft. Unlike simplified models that treat delay as a generic buffer, our framework rigorously calculates the trajectory geometry to ensure that the required separation standards are met at the FAF while maintaining flyable ground speeds and turn radii.

\subsection{Generation of Arrival Flows \label{subsec: flow_generation}}
We simulate arrival traffic originating from four distinct corner fixes surrounding the terminal airspace: DALAS (northwest), LOGEN (northeast), HUSKY (southeast), and TIROE (southwest). Each corner fix generates an independent stream of aircraft modeled by a shifted Poisson process, as in our previous work \cite{pang2025modeling}. This approach incorporates the stochastic nature of arrival demands while strictly enforcing the minimum separation standards required for safety-critical air traffic operations.
\subsubsection{Shifted Poisson Arrival Process}
We denote the set of corner fixes as $\mathcal{K} = \{\text{DALAS}, \text{LOGEN}, \text{HUSKY}, \text{TIROE}\}$. For each corner $\kappa \in \mathcal{K}$, the arrival times are generated as a renewal process where the inter-arrival time $\Delta t$ consists of a deterministic safety buffer and a stochastic component. Let $T_{\mathrm{sep}}$ denote the minimum required separation (66 seconds in our implementation). The inter-arrival time is defined as,
\begin{equation}
\label{eq:interarrival}
\Delta t = T_{\mathrm{sep}} + X_\kappa, \quad X_\kappa \sim \text{Exp}(\lambda_\kappa)
\end{equation}
where $\lambda_\kappa$ is the arrival rate for corner $\kappa$ (in aircraft per hour). The entry time of the $m$-th aircraft at corner $\kappa$, denoted $\tau_{\kappa,m}$, is generated recursively:
\begin{equation}
\label{eq:recursive_entry}
\tau_{\kappa,m} = \tau_{\kappa,m-1} + T_{\mathrm{sep}} + \text{Exp}(\lambda_\kappa), \quad \tau_{\kappa,0} = 0
\end{equation}
This formulation ensures that no two aircraft from the same stream enter the simulation closer than $T_{\mathrm{sep}}$ seconds apart. The generation process terminates when the next computed entry time would exceed a maximum simulation window $T_{\max}$ (3600 seconds in our implementation), ensuring all generated aircraft can be processed within the optimization horizon.
\subsubsection{Traffic Density Sampling}
To evaluate the optimization framework across a wide spectrum of congestion levels, we vary the arrival rate $\lambda_\kappa$ for each corner fix independently in every simulation scenario. Unlike previous approaches that sample rates from a log-uniform distribution, we employ a uniform integer sampling scheme to provide equal probability mass across the specified range of traffic densities.
For each corner $\kappa$, the arrival rate (in aircraft per hour) is sampled as:
\begin{equation}
\label{eq:lambda_sampling}
\lambda_\kappa \sim \mathsf{Uniform}\{\lambda_{\min}, \lambda_{\min}+1, \dots, \lambda_{\max}\}
\end{equation}
where $\lambda_{\min} = 1$ and $\lambda_{\max} = 60$ aircraft per hour in our experiments. This discrete uniform sampling ensures that each integer rate value in the range has equal probability, facilitating straightforward interpretation of results and enabling direct comparison across scenarios with identical demand levels. 

\subsubsection{Global Indexing and Geometry Initialization}
The independent streams from all four corners are aggregated to form the global set of aircraft $\mathcal{I} = \{1, \dots, n\}$ used in the optimization formulation. We sort the total $n$ aircraft by their entry times such that $\tau_1 \le \tau_2 \le \dots \le \tau_n$, where $\tau_i$ denotes the TCP entry time of aircraft $i$.
Crucially, each aircraft $i$ inherits the spatial properties of its generating corner $\kappa$. We define the entry state $C_i = (X_{C_i}, Y_{C_i})$ as the coordinates of the corner fix $\kappa$ where aircraft $i$ originated. This entry position $C_i$ determines the specific geometric flight path characteristics utilized in the travel time calculation in Section \ref{subsec: travel_time}.

\subsection{TCP Travel Time \label{subsec: travel_time}}
We define the aircraft $i$'s location when entering the TCP at time $\tau_i$ as $C_i = (X_{C_i}, Y_{C_i})$. The target is the Final Approach Fix (FAF) at $(X_{FAF}, Y_{FAF})$. The flight path consists of a tangent leg $d_L$, an RF turn arc $d_\theta$ with radius $r$, and a final straight-in offset $d_i$, which serves as one of our decision variables ($d_i \ge 0$).

The geometry of the path depends on the aircraft's position relative to the FAF when entering the TCP. The center of the RF turn, $C_{0, i}$, and the reference projection $C'_{0,i}$ are defined based on the decision variable $d_i$ and the arrival flow direction (either north arrivals or south arrivals),
\begin{equation}
C_{0,i}(d_i)=
\begin{cases}
\bigl(X_{FAF}-d_i,\; Y_{FAF}+r\bigr), & \text{if } Y_{C_i}>Y_{FAF}+r,\\
\bigl(X_{FAF}-d_i,\; Y_{FAF}-r\bigr), & \text{if } Y_{C_i}<Y_{FAF}-r.
\end{cases}
\end{equation}

\begin{equation}
C'_{0,i}(d_i)=\bigl(X_{FAF}-d_i,\; Y_{FAF}\bigr).
\end{equation}

The total flight path length $D_i(d_i)$ is the summation $D_i(d_i) = d_L(d_i) + d_\theta(d_i) + d_i$. For $d_L$, further defining the Euclidean norm $d_0 = ||C_i-C_{0,i}(d_i)||_2$, we have $d_L = \sqrt{d_0^2-r^2}$. For the arc length $d_\theta$, we need to compute the corresponding angle $\theta$ of the arc. Further define $\theta_1 = \arccos (\frac{r^2+d_0^2-d_L^2}{2 r d_0})$, and $\theta_2 = \arccos (\frac{r^2+d_0^2-d_0^{'2}}{2 r d_0})$, where $d_0'$ is the Euclidean norm as $d_0' = ||C_i - C'_{0,i}(d_i)||_2$. For all east arrivals, $\theta = 2\pi - (\theta_1 + \theta_2)$, while $\theta = \theta_2 - \theta_1$ for all west arrivals. $d_\theta$ is simply the fraction of the arc as $d_\theta = \theta r$. Giving a commanded constant speed $v_L$ for the tangent leg $d_L$, a constant ground speed $v_\theta$ for the RF leg, and speed $v_i$ for the straight-in leg, the total travel time $t_i$ is given by,
\begin{equation}
\label{eq: travel_time}
   t_i = \tau_i + \frac{d_L(d_i)}{v_L} + \frac{d_\theta(d_i)}{v_\theta} + \frac{d_i}{v_i}
\end{equation}
\noindent where $\tau_i$ is the initial entry time.

It is worth pointing out that the derived quantities $d_0$, $d'_0$, $d_L$, $\theta_1$, $\theta_2$, and $d_\theta$ all depend on $d_i$. The complete TCP travel time formulation for each arrival aircraft $i$ is derived as follows,

In this formulation, we assume the following parameters,
\begin{itemize}
    \item $r$ is the RF turn radius of the arrival aircraft. In practice, this varies for different aircraft types at different airports. 
    \item $X_{FAF}$ and $Y_{FAF}$ are the coordinates of the FAF of the landing runway. In our setup, we use $\mathsf{BURNY}$ as the FAF for runway 9R at KATL.
\end{itemize}

\begin{figure}[H]
\centering
    \begin{subfigure}[t]{0.45\textwidth}
        \includegraphics[width=\textwidth]{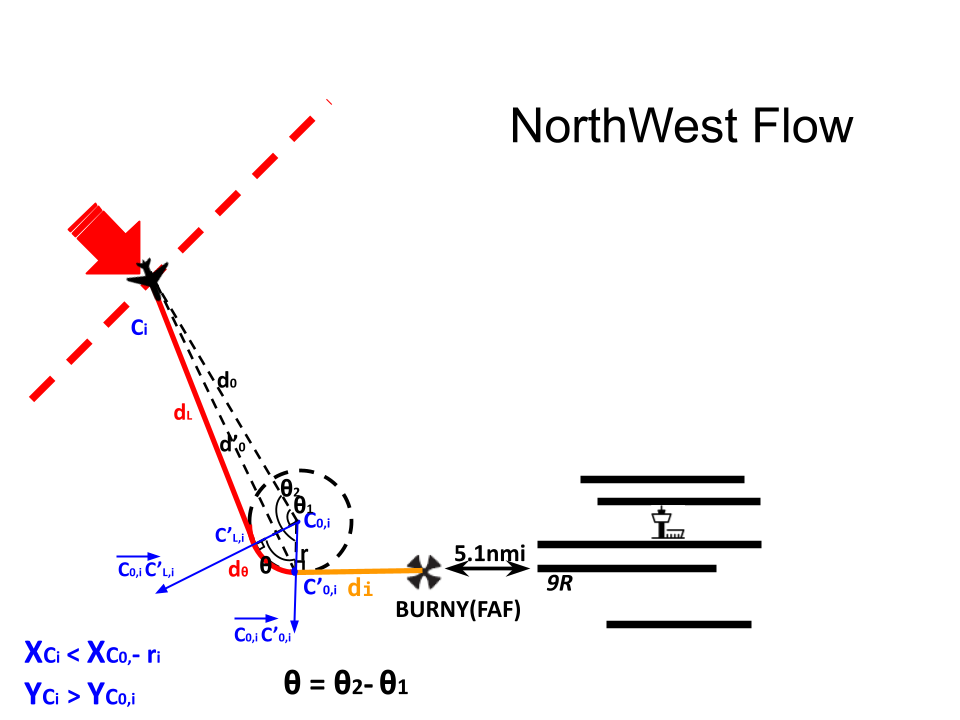}
        \caption{Landing geometry of northwest arrivals.}
        \label{fig: NorthWestFlow}
    \end{subfigure}
    \begin{subfigure}[t]{0.45\textwidth}
        \includegraphics[width=\textwidth]{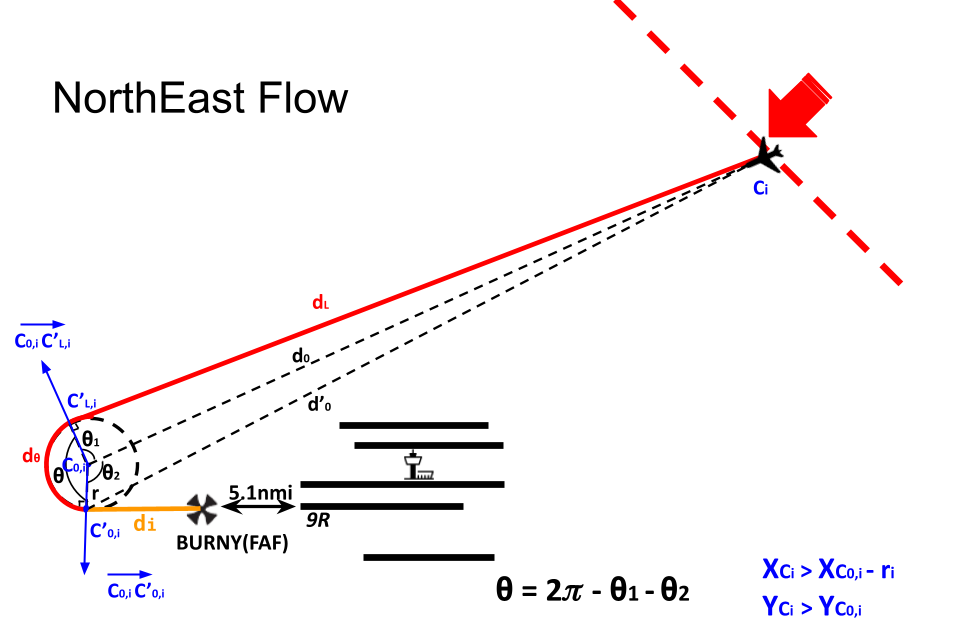}
        \caption{Landing geometry of northeast arrivals.}
        \label{fig: NorthWestFlow}
    \end{subfigure}
    \\
    \begin{subfigure}[t]{0.45\textwidth}
        \includegraphics[width=\textwidth]{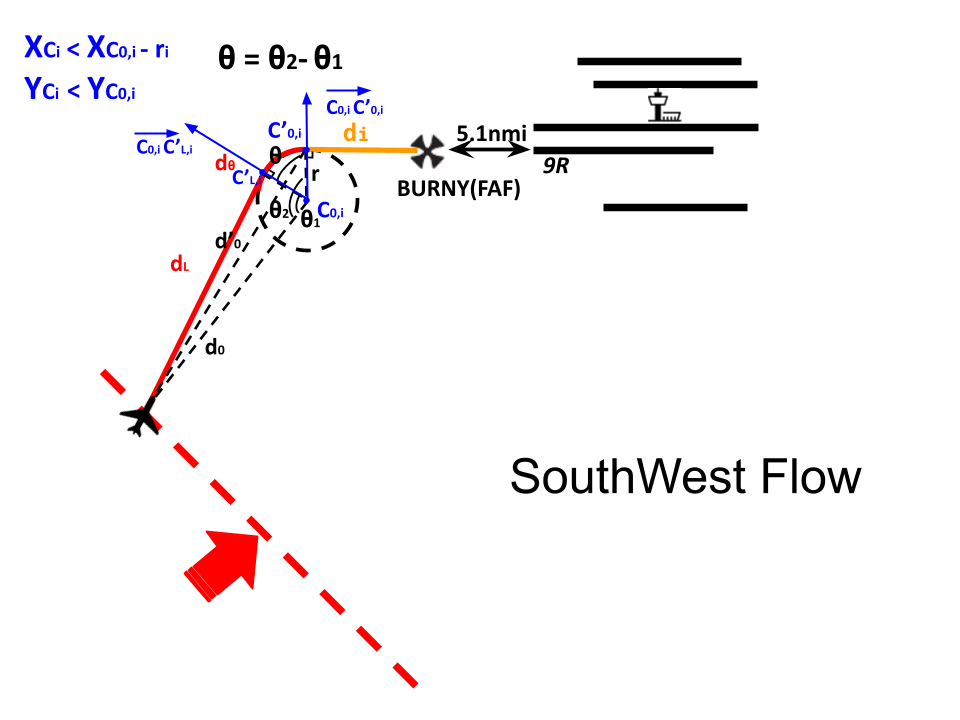}
        \caption{Landing geometry of southwest arrivals.}
        \label{fig: SouthWestFlow}
    \end{subfigure}
    \begin{subfigure}[t]{0.45\textwidth}
        \includegraphics[width=\textwidth]{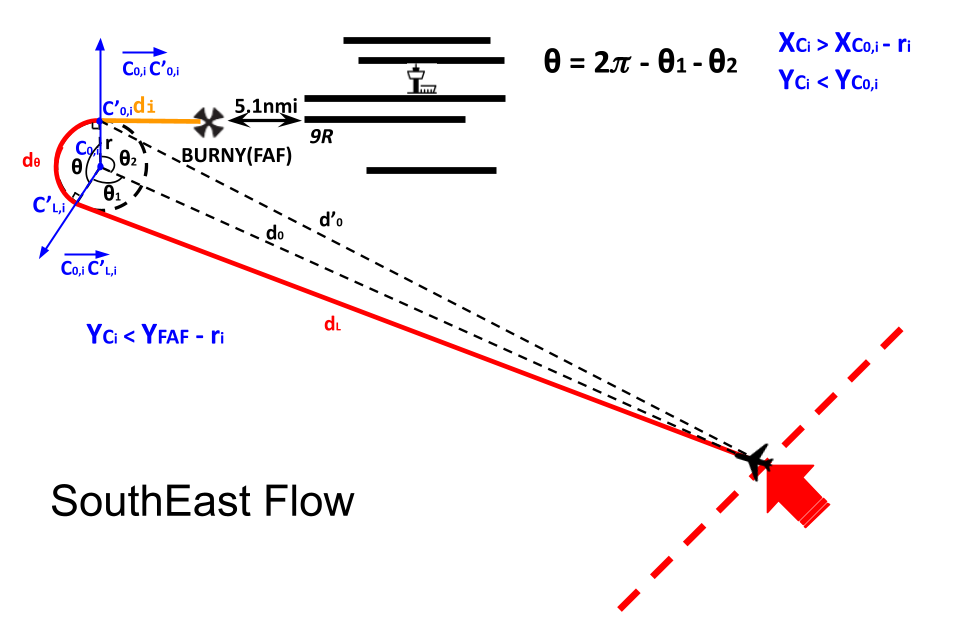}
        \caption{Landing geometry of southeast arrivals.}
        \label{fig: SouthWestFlow}
    \end{subfigure}
\caption[Angle calculation using vectors for the four arrival flows.]{We can calculate $\theta_2(d_i)$ using the angle between vectors. When $Y_{C_i}<Y_{C_{0,i}}$, we have $\theta_2$ to be the clock wise angle from vector $\overrightarrow{C_{0,i}C_i} = (X_{FAF}-d_i-X_{C_i}, Y_{FAF}-r-Y_{C_i})$ to vector $\overrightarrow{C_{0,i}C'_{0,i}} = (0, r)$, and when $Y_{C_i}>Y_{C_{0,i}}$, we have counter clockwise angle from vector $\overrightarrow{C_{0,i}C_i} = (X_{FAF}-d_i-X_{C_i}, Y_{FAF}+r-Y_{C_i})$ to vector $\overrightarrow{C_{0,i}C'_{0,i}} = (0, -r)$.}
\label{fig: geometry}
\end{figure}

\subsection{Computation of RF Arch Segment Length $d_{\theta_i(d_i)}$}
We can directly compute the RF-leg central angle $\theta(d_i)$ using the angle between two vectors, without explicitly evaluating $\theta_2(d_i)$ via angle wrapping. The key step is to obtain the coordinate of the left tangent point, denoted by $C_{L,i}(d_i) = \bigl(X_{C_{L,i}}(d_i),\, Y_{C_{L,i}}(d_i)\bigr)$, which is the unique tangency point on the RF circle such that the tangent segment $\overrightarrow{C_i C_{L,i}}$ is tangent to the circle and $X_{C_{L,i}}$ is minimized (i.e., the left tangent). That is, we have,
\begin{equation}
    \overrightarrow{C_i C_{L,i}} \cdot \overrightarrow{ C_{L,i}C_{0,i}}=0
\end{equation}
which implies that the line segment $\overrightarrow{C_i C_{L,i}}$ is perpendicular to the radius $\overrightarrow{C_{L,i} C_{0,i}}$ at the tangent point.

Define the center-to-point vector,
\begin{equation}
\label{eq: center_to_point_vector}
\mathbf v_i(d_i)=C_i-C_{0,i}(d_i)=
\begin{bmatrix}
X_{C_i}-X_{C_{0,i}}(d_i)\\
Y_{C_i}-Y_{C_{0,i}}(d_i)
\end{bmatrix},
\end{equation}
Note that $\mathbf v_i(d_i)$ is \emph{not} a tangent direction; it is used only as a construction basis. A valid tangent exists only if $d_0(d_i)>r$. Define the perpendicular basis vector
\begin{equation}
\label{eq: perp_basis_vector}
\mathbf v_i^\perp(d_i)=
\begin{bmatrix}
-\bigl(Y_{C_i}-Y_{C_{0,i}}(d_i)\bigr)\\
\phantom{-}\bigl(X_{C_i}-X_{C_{0,i}}(d_i)\bigr)
\end{bmatrix},
\end{equation}
and introduce the coefficients
\begin{equation}
\label{eq: tangent_coefficients}
a_i(d_i)=\frac{r^2}{d_0(d_i)^2},
\qquad
b_i(d_i)=\frac{r\sqrt{d_0(d_i)^2-r^2}}{d_0(d_i)^2}.
\end{equation}
Then the left tangent point $C_{L,i}(d_i)$ can be written in closed form as
\begin{equation}
\label{eq: left_tangent_point_two_cases}
C_{L,i}(d_i)=
\begin{cases}
C_{0,i}(d_i)+a_i(d_i)\,\mathbf v_i(d_i)+b_i(d_i)\,\mathbf v_i^\perp(d_i), & \text{if } Y_{C_i} > Y_{FAF} + r, \\
C_{0,i}(d_i)+a_i(d_i)\,\mathbf v_i(d_i)-b_i(d_i)\,\mathbf v_i^\perp(d_i), & \text{if } Y_{C_i} < Y_{FAF} - r.
\end{cases}
\end{equation}
Equivalently, the coordinates of the tangent point $C_{L,i}(d_i)=\bigl(X_{C_{L,i}}(d_i),Y_{C_{L,i}}(d_i)\bigr)$ satisfy
\begin{equation}
\label{eq: left_tangent_coordinates_two_cases}
C_{L,i}(d_i)=
\begin{cases}
\Bigl(X_{C_{0,i}}(d_i)+a_i\Delta X_i-b_i\Delta Y_i,\;\; Y_{C_{0,i}}(d_i)+a_i\Delta Y_i+b_i\Delta X_i\Bigr), & \text{if } Y_{C_i} > Y_{FAF} + r \\
\Bigl(X_{C_{0,i}}(d_i)+a_i\Delta X_i+b_i\Delta Y_i,\;\; Y_{C_{0,i}}(d_i)+a_i\Delta Y_i-b_i\Delta X_i\Bigr), & \text{if } Y_{C_i} < Y_{FAF} - r \\
\end{cases}
\end{equation}
where $\Delta X_i=X_{C_i}-X_{C_{0,i}}(d_i)$ and $\Delta Y_i=Y_{C_i}-Y_{C_{0,i}}(d_i)$, and we use the shorthand $a_i=a_i(d_i)$ and $b_i=b_i(d_i)$ for compactness. 

With the left tangent point available, the RF-leg central angle $\theta_i(d_i)$ can be computed directly from the angle between the two radius vectors as,
\begin{equation}
\label{eq: theta_via_dot_cross}
\theta_i(d_i)=\operatorname{atan2}\!\left(
\left|\overrightarrow{C_{0,i}C'_{0,i}} \times \overrightarrow{C_{0,i}C_{L,i}}\right|,
\overrightarrow{C_{0,i}C'_{0,i}} \cdot \overrightarrow{C_{0,i}C_{L,i}}
\right)\in[0,\pi].
\end{equation}

Using this equation and the fact that $\overrightarrow{C_{0,i}C'_{0,i}}$ is purely vertical, the RF-leg central angle can be written directly in terms of the left tangent point coordinates. Specifically,

\begin{equation}
\label{eq: theta_final_no_deltas_using_CL}
\theta_i(d_i)=
\begin{cases}
\operatorname{atan2}\!\Bigl(\bigl|X_{C_{L,i}}(d_i)-X_{C_{0,i}}(d_i)\bigr|,\; -\bigl(Y_{C_{L,i}}(d_i)-Y_{C_{0,i}}(d_i)\bigr)\Bigr),
& \text{if } Y_{C_i} > Y_{\mathrm{FAF}}+r,\\[6pt]
\operatorname{atan2}\!\Bigl(\bigl|X_{C_{L,i}}(d_i)-X_{C_{0,i}}(d_i)\bigr|,\;\;\;\bigl(Y_{C_{L,i}}(d_i)-Y_{C_{0,i}}(d_i)\bigr)\Bigr),
& \text{if } Y_{C_i} < Y_{\mathrm{FAF}}-r,
\end{cases}
\end{equation} 

which yields $\theta_i(d_i)\in[0,\pi]$ by construction. The corresponding RF arc length is $d_{\theta,i}(d_i)=r\,\theta_i(d_i)$.

Substituting the closed-form expression of $C_{L,i}(d_i)$ from \eqref{eq: left_tangent_coordinates_two_cases} further gives $\theta_i(d_i)$ explicitly in terms of $(X_{C_i},Y_{C_i})$, $(X_{\mathrm{FAF}},Y_{\mathrm{FAF}})$, $r$, and $d_i$. Let
\[
X_{C_{0,i}}(d_i)=X_{\mathrm{FAF}}-d_i,\qquad
Y_{C_{0,i}}(d_i)=
\begin{cases}
Y_{\mathrm{FAF}}+r, & Y_{C_i} > Y_{\mathrm{FAF}}+r,\\
Y_{\mathrm{FAF}}-r, & Y_{C_i} < Y_{\mathrm{FAF}}-r,
\end{cases}
\]
and define
\[
d_0(d_i)=\sqrt{\bigl(X_{C_i}-X_{C_{0,i}}(d_i)\bigr)^2+\bigl(Y_{C_i}-Y_{C_{0,i}}(d_i)\bigr)^2},\qquad d_0(d_i)>r,
\]
\[
a_i(d_i)=\frac{r^2}{d_0(d_i)^2},\qquad
b_i(d_i)=\frac{r\sqrt{d_0(d_i)^2-r^2}}{d_0(d_i)^2}.
\]
Then,
\begin{equation}
\label{eq: theta_final_no_deltas_expanded}
\theta_i(d_i)=
\begin{cases}
\operatorname{atan2}\!\Bigl(
\bigl|\,a_i(d_i)\bigl(X_{C_i}-X_{\mathrm{FAF}}+d_i\bigr)-b_i(d_i)\bigl(Y_{C_i}-Y_{\mathrm{FAF}}-r\bigr)\,\bigr|, \\
\;-\Bigl[a_i(d_i)\bigl(Y_{C_i}-Y_{\mathrm{FAF}}-r\bigr)+b_i(d_i)\bigl(X_{C_i}-X_{\mathrm{FAF}}+d_i\bigr)\Bigr]
\Bigr),
& \text{if } Y_{C_i} > Y_{\mathrm{FAF}}+r,\\[8pt]
\operatorname{atan2}\!\Bigl(
\bigl|\,a_i(d_i)\bigl(X_{C_i}-X_{\mathrm{FAF}}+d_i\bigr)+b_i(d_i)\bigl(Y_{C_i}-Y_{\mathrm{FAF}}+r\bigr)\,\bigr|, \\
\;\;\;\Bigl[a_i(d_i)\bigl(Y_{C_i}-Y_{\mathrm{FAF}}+r\bigr)-b_i(d_i)\bigl(X_{C_i}-X_{\mathrm{FAF}}+d_i\bigr)\Bigr]
\Bigr),
& \text{if } Y_{C_i} < Y_{\mathrm{FAF}}-r,
\end{cases}
\end{equation}

The corresponding RF arc length is given by $d_{\theta_i(d_i)}=r\,\theta(d_i)$.

\subsection{Fixed-Sequence Trajectory Optimization}
\label{subsec:baseline_nlp}
We consider a set of arrival aircraft $\mathcal{I}=\{1,\dots,N\}$ and a fixed landing sequence represented by the permutation vector $\Pi,$ where $\Pi(k)$ denotes the aircraft index scheduled to land in the $k$-th position. The sequence $\Pi$ is determined a priori using FCFS ordering based on nominal ETAs and is treated as a hard constraint throughout the optimization. The decision variables are each aircraft's path extension $d_i$ and segment-wise speed profiles $(v_{L,i}, v_{\theta,i}, v_{f,i})$ to satisfy separation requirements while minimizing system-level objectives.

\subsubsection{Decision variables}
For each $i\in\mathcal{I},$ we optimize
\begin{equation}
\label{eq:vars_nlp}
d_i \in [0,D_{\max}],\
v_{L,i}\in[V_{L}^{\min},V_{L}^{\max}], \
v_{\theta,i}\in[V_{\theta}^{\min},V_{\theta}^{\max}],\
v_{f,i}\in[V_{f}^{\min},V_{f}^{\max}],\
t_i\ge 0,
\end{equation}
and for each rank $k=2,\dots,N,$ we introduce a soft separation slack variable $\sigma_k \ge 0.$.

\subsubsection{Constraints} 
\noindent\textbf{(i) Travel-time dynamics (FAF arrival time).}
The arrival time $t_i$ at the FAF is determined by the path geometry and segment speeds:
\begin{equation}
\label{eq:nlp_time}
t_i
\geq
\tau_i
+
\frac{d_L(d_i)}{v_{L,i}}
+
\frac{r\,\theta(d_i)}{v_{\theta,i}}
+
\frac{d_i}{v_{f,i}},
\qquad \forall i\in\mathcal{I},
\end{equation}
where $d_L(d_i)$ and $\theta(d_i)$ are computed using the closed-form expressions from Section \ref{subsec: travel_time} that account for the trombone path extension $d_i$.

\textbf{(ii) Fixed-sequence separation (soft).}
Let $t_k := t_{\Pi(k)}$ denote the FAF arrival time of the aircraft at rank $k$. The separation constraint is enforced as a soft constraint:
\begin{equation}
\label{eq:nlp_sep}
t_k \ge t_{k-1}+T_{\mathrm{sep}}-\sigma_k,
\qquad \forall k\in\{2,\dots,N\},
\end{equation}
where $\sigma_k \ge 0$ is a slack variable that allows violations but incurs a large penalty $M_{\mathrm{pen}}$ in the objective.

\textbf{(iii) Speed monotonicity across segments.}
To ensure realistic flight profiles, speeds must decrease monotonically from the tangent leg to the final approach:
\begin{equation}
\label{eq:nlp_speedmono}
v_{L,i} \ge v_{\theta,i} \ge v_{f,i},
\qquad \forall i\in\mathcal{I}.
\end{equation}
\noindent\textbf{(iv) Variable bounds.}
\begin{equation}
\label{eq:nlp_bounds}
0\le d_i \le D_{\max},
\
V_{\ell}^{\min}\le v_{\ell,i}\le V_{\ell}^{\max}\ \ (\ell\in\{L,\theta,f\}),
\
t_i\ge 0,
\
\sigma_k\ge 0.
\end{equation}
\subsubsection{Objective Function}
\label{sec:single_stage_obj}

We propose a unified single-stage optimization approach that simultaneously addresses safety, throughput, and efficiency. We formulate a scalar objective function $J^\star$ as a weighted sum of three competing goals. This approach maintains the operational hierarchy of Air Traffic Control by assigning dominant weights to safety constraints, followed by throughput, and finally efficiency.

The global optimization problem is defined as,

\begin{equation}
\label{eq:single_stage_objective}
J^\star
=
\min_{\substack{\{d_i,v_{L,i},v_{\theta,i},v_{f,i},t_i\} \\ \{\sigma_k\}}}
\left(
W_{\mathrm{safe}} \sum_{k=2}^{N} \sigma_k 
\;+\; W_{\mathrm{thru}} \cdot t_N 
\;+\; W_{\mathrm{eff}} \sum_{i=1}^N d_i
+ W_{\mathrm{speed}}\sum_{i\in\mathcal I}\sum_{\ell\in\{L,\theta,f\}}
\frac{V_\ell^{\max}-v_{\ell,i}}{V_\ell^{\max}-V_\ell^{\min}}.
\right)
\end{equation}

\noindent\textbf{Subject to:}
\begin{itemize}
    \item Constraints \eqref{eq:nlp_time}--\eqref{eq:nlp_bounds}
\end{itemize}

\noindent\textbf{Where:}
\begin{itemize}
    \item $W_{\mathrm{safe}}$: A dominant penalty weight ensuring separation violations $\sigma_k$ are driven to zero whenever a feasible solution exists.
    \item $W_{\mathrm{thru}}$: A primary operational weight to minimize the makespan $t_N$ (time of the last FAF crossing).
    \item $W_{\mathrm{eff}}$: A secondary weight to minimize total path stretch $\sum d_i$, prioritizing fuel efficiency only after safety and throughput targets are met.
    \item $W_{\mathrm{speed}}$: A small weight to prefer higher speed for fuel efficiency, biases toward upper speed bounds when other objectives are satisfied. 
    \item $\sigma_k \ge 0$: Slack variables representing the magnitude of separation violation for the $k$-th aircraft.
\end{itemize}

\section{Simulations and Results}\label{sec: results}
\subsection{Case Study with A80 TRACON}
We evaluate the proposed single-stage optimization framework through comprehensive Monte Carlo simulations using the Atlanta A80 TRACON airspace structure. The simulation environment models arrival traffic from four feeder gates (DALAS, LOGEN, HUSKY, TIROE) using the shifted Poisson process described in Section~\ref{subsec: flow_generation}. The simulation parameters are configured to reflect realistic terminal operations. The simulation time horizon is one hour, and only aircraft with entry times within this window are processed to ensure computational tractability.

\begin{figure}
\centering
    \begin{subfigure}[t]{0.45\textwidth}
        \includegraphics[width=\textwidth]{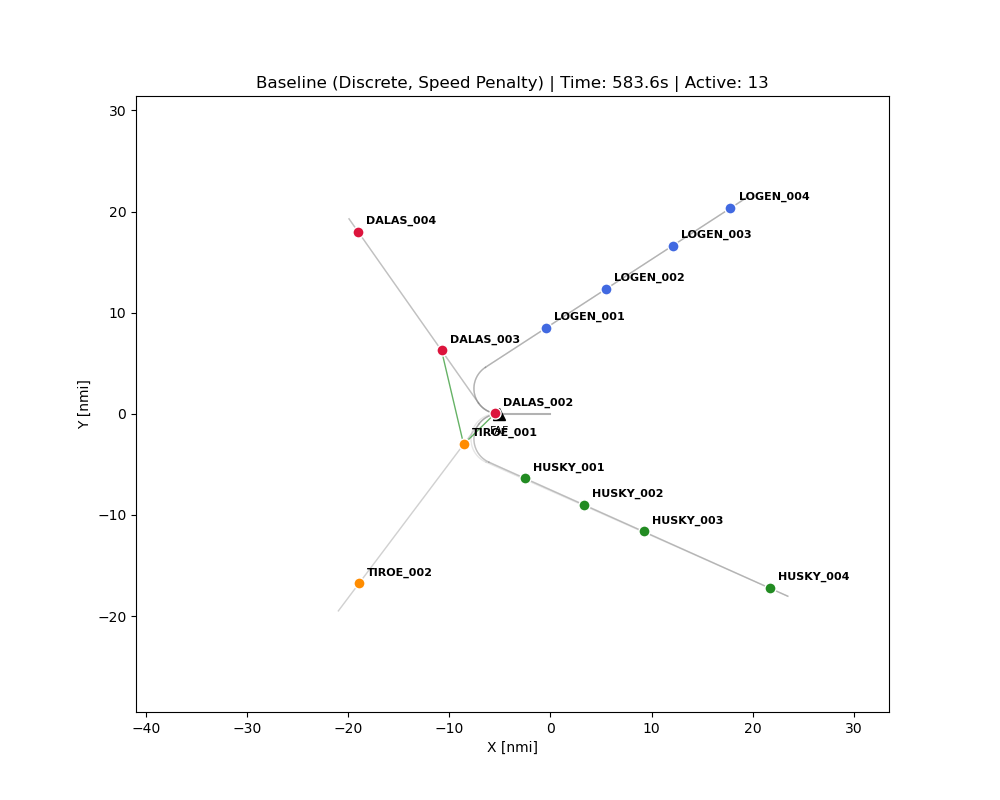}
        \caption{FCFS with discrete speed decreasing penalty at simulation time $t=583.6$s.}
        \label{fig: fcfs-frame1}
    \end{subfigure}
    \begin{subfigure}[t]{0.45\textwidth}
        \includegraphics[width=\textwidth]{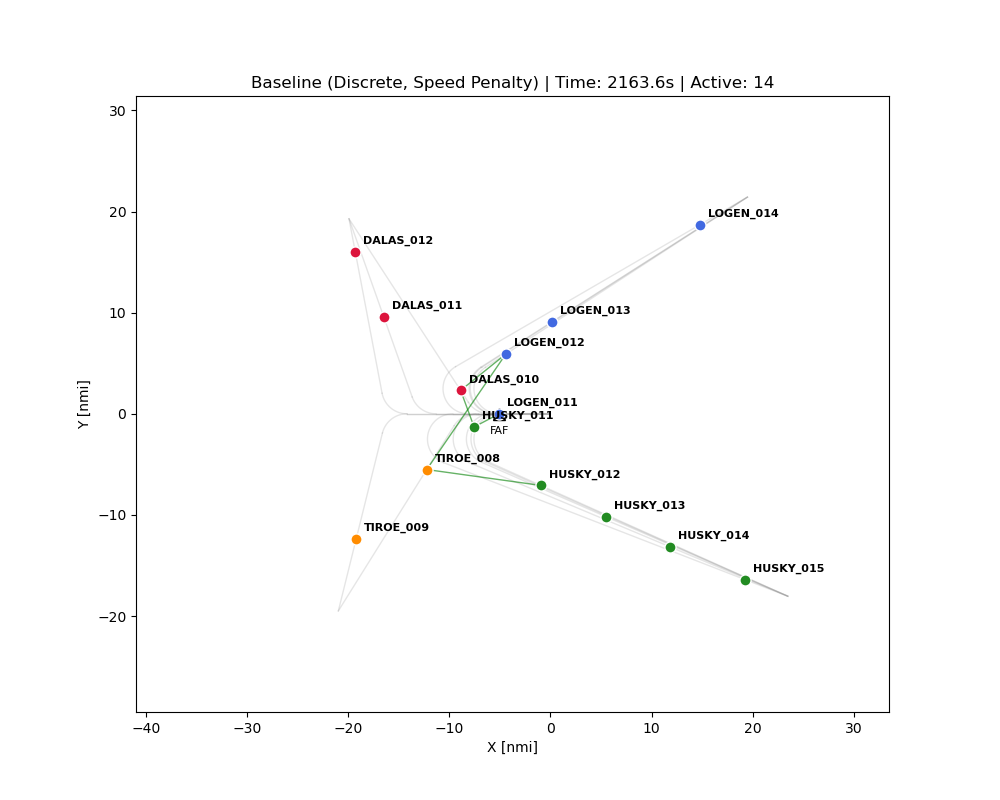}
        \caption{FCFS with discrete speed decreasing penalty at simulation time $t=2163.6$s.}
        \label{fig: fcfs-frame2}
    \end{subfigure}
    \\
    \begin{subfigure}[t]{0.45\textwidth}
        \includegraphics[width=\textwidth]{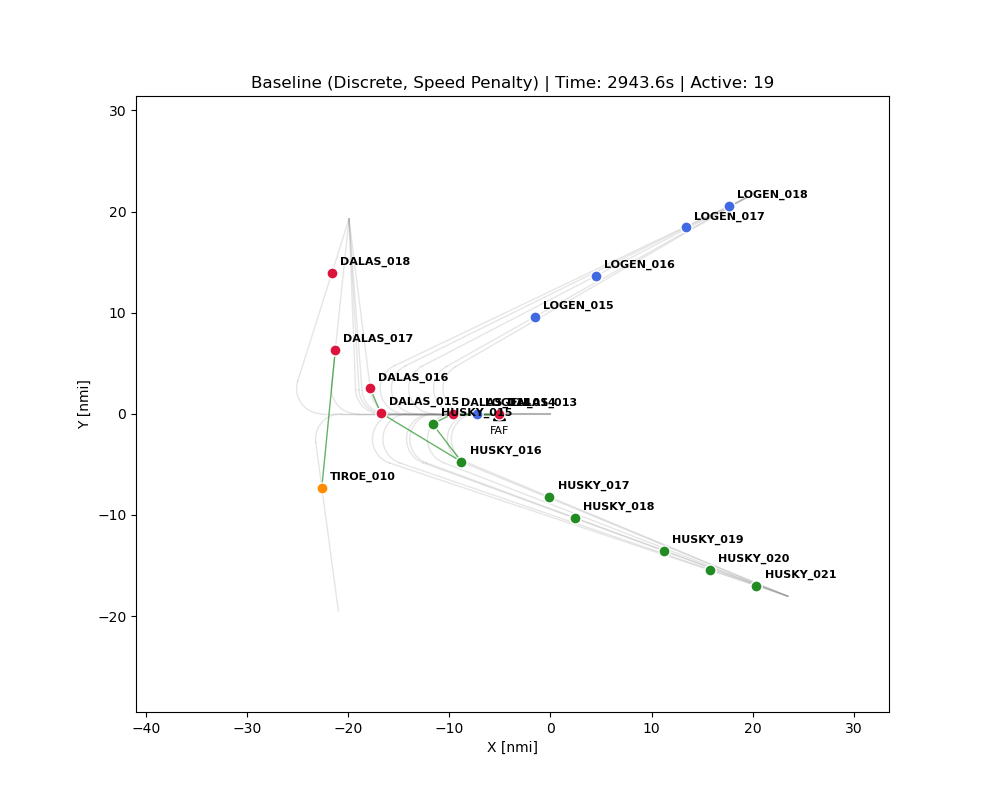}
        \caption{FCFS with discrete speed decreasing penalty at simulation time $t=2943.6$s.}
        \label{fig: fcfs-frame3}
    \end{subfigure}
    \begin{subfigure}[t]{0.45\textwidth}
        \includegraphics[width=\textwidth]{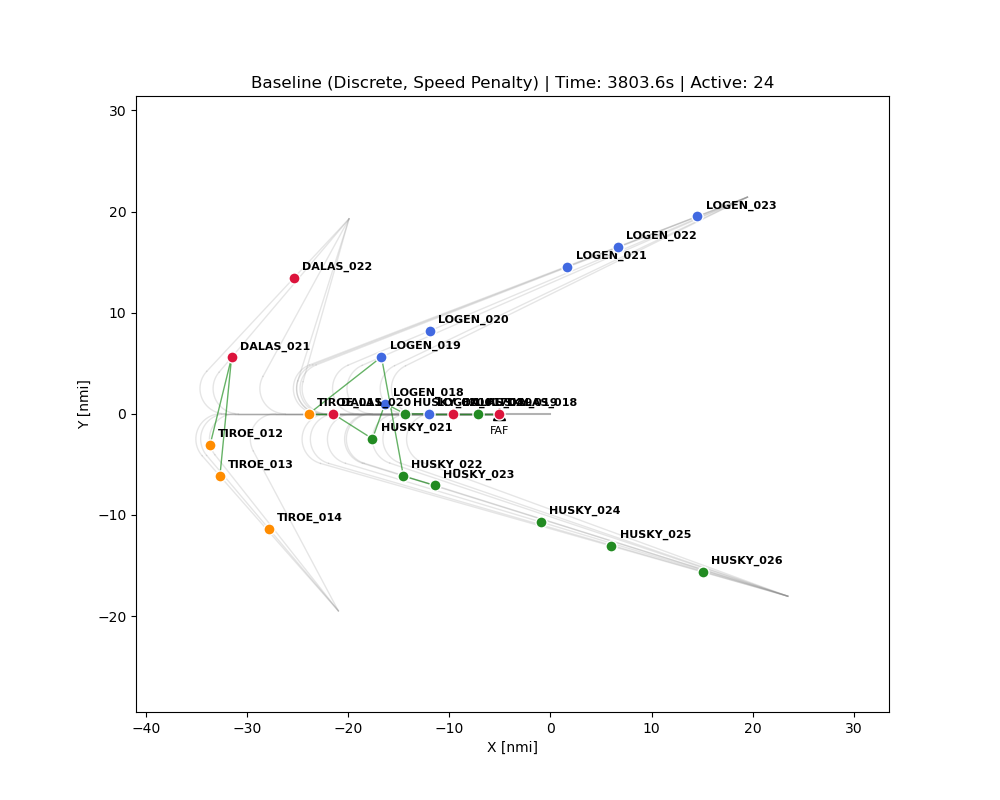}
        \caption{FCFS with discrete speed decreasing penalty at simulation time $t=3803.6$s.}
        \label{fig: fcfs-frame4}
    \end{subfigure}
    \\
    \begin{subfigure}[t]{0.45\textwidth}
        \includegraphics[width=\textwidth]{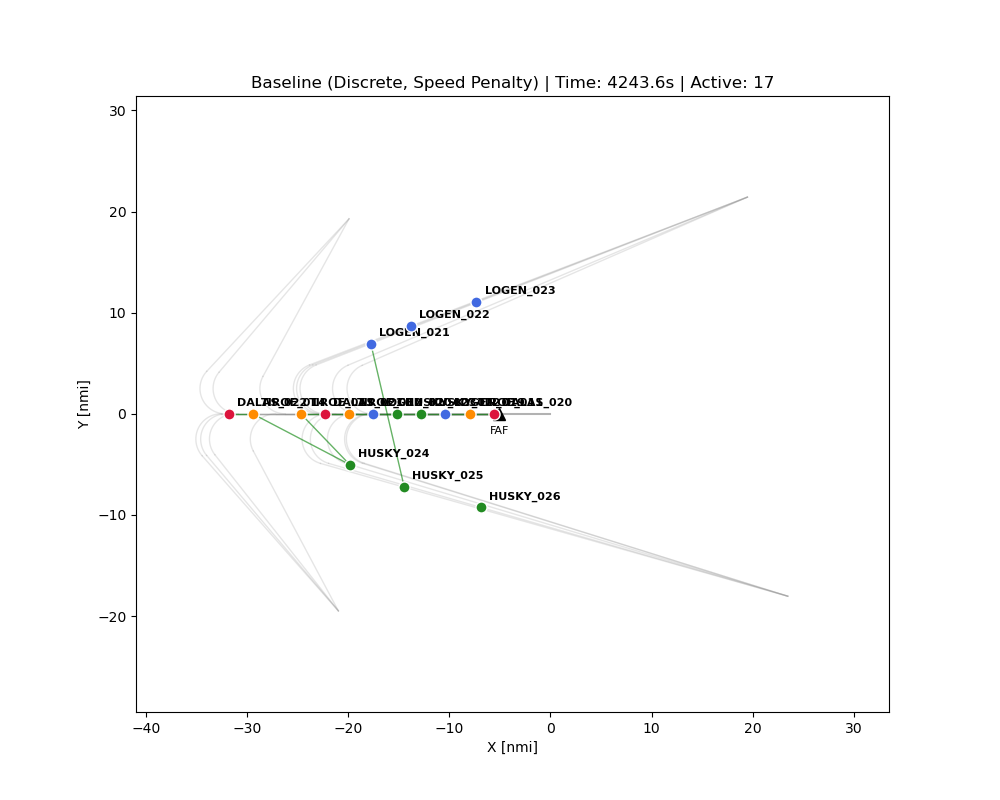}
        \caption{FCFS with discrete speed decreasing penalty at simulation time $t=4243.6$s.}
        \label{fig: fcfs-frame5}
    \end{subfigure}
    \begin{subfigure}[t]{0.45\textwidth}
        \includegraphics[width=\textwidth]{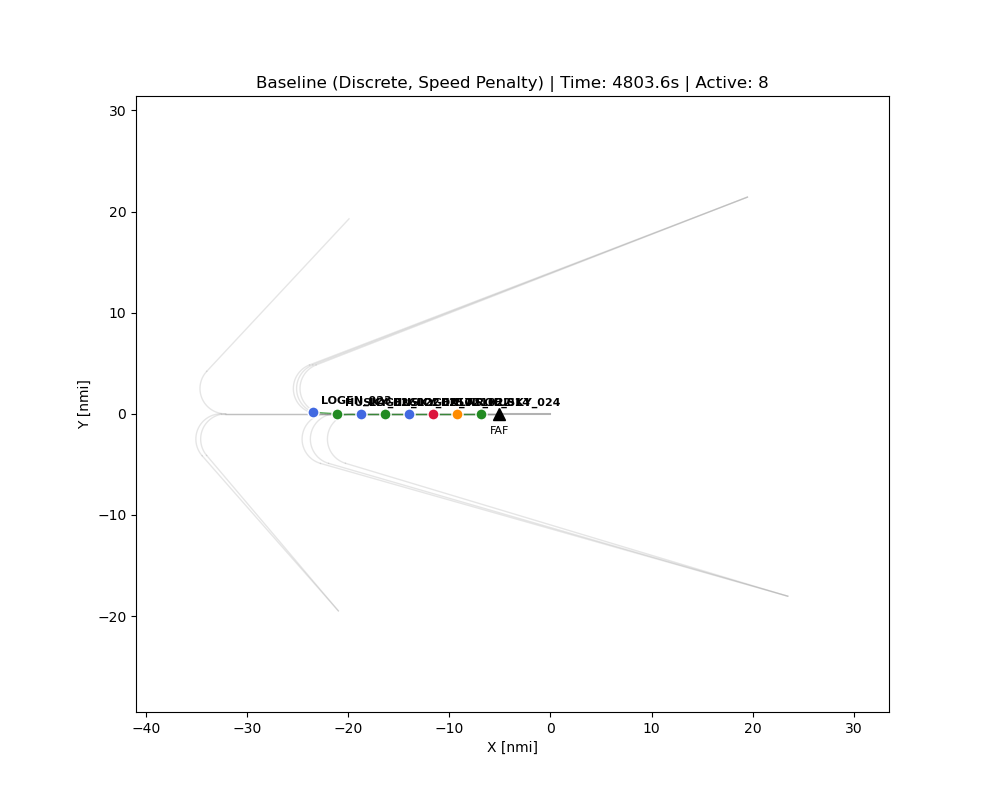}
        \caption{FCFS with discrete speed decreasing penalty at simulation time $t=4803.6$s.}
        \label{fig: fcfs-frame6}
    \end{subfigure}
\caption{Visualization of the optimization results for one random generated scenario. Green line connects the successive landings that are well-separated. In this scenario, the speeds take 10-knot intervals, with penalty on decreased speed from the upper bound, simulating the first-come-first-serve (FCFS) logic for terminal arrivals. It is clear that in order to create the required separation at the Baseleg extensions without separation violations, the $d_i$s are pushed to the maximum value. Runway threshold is the origin.}
\label{fig: fcfs-vis}
\end{figure}

Speed bounds are defined to reflect realistic approach profiles. For the tangent leg segment, speeds are constrained to $v_{L,i} \in [180, 240]$ knots. For the RF turn segment, speeds range from $v_{\theta,i} \in [130, 200]$ knots, allowing for reduced speeds during the curved approach. The final approach segment uses speeds $v_{f,i} \in [130, 160]$ knots, consistent with stabilized approach requirements. These bounds ensure that the generated trajectories maintain realistic flight profiles while providing sufficient control authority for delay absorption.

The optimization is solved using IPOPT \cite{biegler2009large} in Pyomo with the weighted objective function defined in Equation~\eqref{eq:single_stage_objective}. The weight hierarchy is configured as $W_{\mathrm{safe}} = 10^4$, $W_{\mathrm{thru}} = 1.0$, $W_{\mathrm{eff}} = 0.1$, and $W_{\mathrm{speed}} = 0.01$ to enforce the operational priority of safety over throughput and efficiency. This weight selection ensures that separation violations are driven to zero whenever physically feasible, while throughput (makespan minimization) takes precedence over fuel efficiency (path stretch minimization) when safety constraints are satisfied.

\Cref{fig: fcfs-vis} visualizes the progression of a high-density arrival scenario, showing six snapshots of aircraft positions and trajectories at different simulation times. The visualization demonstrates the algorithm's behavior under congestion, where aircraft are forced to extend their Baseleg paths to maintain required separation standards. 

\subsection{Monte Carlo Tests}

\begin{figure}[t]
    \centering
    \includegraphics[width=0.75\textwidth]{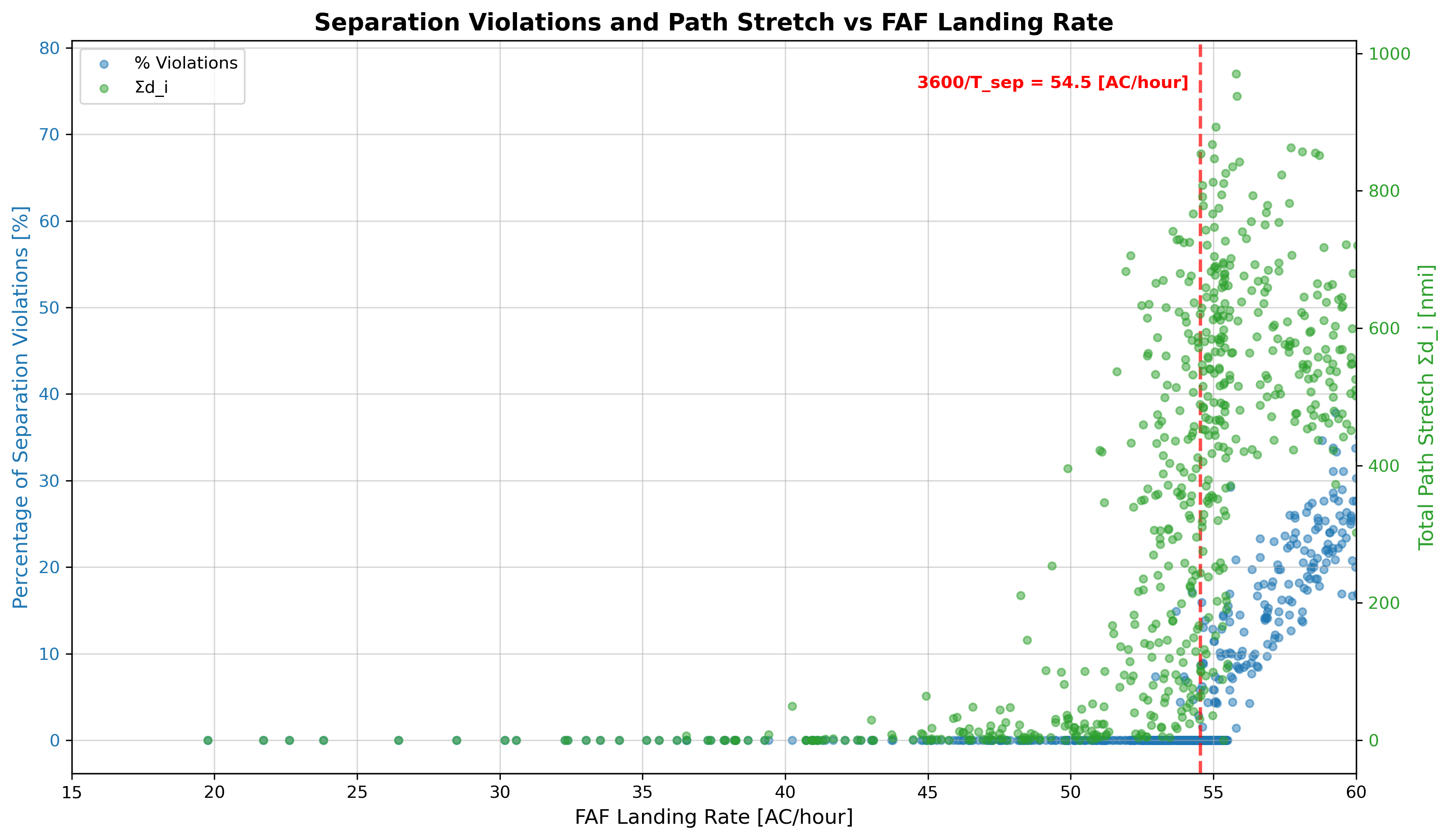}
    \caption{Monte Carlo test of the algorithm with $1,000$ runs. When the density is increased, the algorithm first increases the total stretch to create sufficient separation between landing aircraft. When the arrival rate is higher than the principal maximum tolerable arrival rate at the runway threshold ($3600/T_{\mathrm{sep}}$), the separation violation slack variable increases, leading to separation violations.}
    \label{fig: analysis}
\end{figure}

To quantify the algorithm's performance across a wide spectrum of traffic densities, we conducted a Monte Carlo test with $1,000$ independent simulation runs, each with randomly sampled arrival rates. \Cref{fig: analysis} presents the relationship between the FAF landing rate (the effective throughput achieved at the runway threshold) and two key performance metrics: the percentage of separation violations and the total path stretch $\sum d_i$.

The results reveal a clear two-phase behavior of the optimization algorithm. In the first phase, when the arrival rate is below the theoretical maximum tolerable rate at the runway threshold ($3600/T_{\mathrm{sep}} \approx 54.5$ arrivals per hour), the algorithm successfully maintains zero separation violations by strategically increasing the total path stretch. As shown in \Cref{fig: analysis}, the total stretch $\sum d_i$ increases monotonically with the FAF landing rate, demonstrating the algorithm's ability to absorb delay through geometric path extensions. This behavior aligns with operational ATC practices, where controllers first attempt to manage spacing through speed control, then resort to vectoring (path stretching) when speed adjustments are insufficient.

However, when the arrival rate exceeds the maximum tolerable rate, the system becomes saturated, and the separation violation slack variables $\sigma_k$ begin to increase. As illustrated in \Cref{fig: analysis}, the percentage of separation violations remains near zero for arrival rates below the threshold, then increases sharply once the system capacity is exceeded. This transition point is clearly marked by the vertical dashed line in the figure, representing the fundamental physical limit imposed by the minimum separation requirement $T_{\mathrm{sep}} = 66$ seconds. Beyond this threshold, no amount of path stretching or speed control can maintain separation for all aircraft pairs, and the soft constraint mechanism allows the optimization to find the best achievable solution while explicitly quantifying the magnitude of violations through the slack variables.

The Monte Carlo analysis demonstrates that the single-stage weighted objective function effectively enforces the operational hierarchy: safety (i.e., separation violations) is maintained whenever physically possible, throughput is maximized subject to safety constraints, and efficiency (i.e., minimized path stretch) is optimized only after safety and throughput objectives are satisfied. The smooth transition in the percentage of violations near the capacity threshold indicates that the weight hierarchy $W_{\mathrm{safe}} \gg W_{\mathrm{thru}} \gg W_{\mathrm{eff}}$ successfully approximates lexicographic optimization while maintaining computational efficiency through a single solver call.

\section{Conclusion}\label{sec: conclusion}
This paper presented a high-fidelity trajectory-based optimization framework for the TAAS problem that captures the geometric reality of controller vectoring in high-density TRACON operations (i.e., FCFS based on EAT). Unlike node-link scheduling models that reduce trajectories to time-delay abstractions, we explicitly model the coupling between Baseleg extension ($d_i$) and RF-turn geometry. 

Simulations based on the KATL A80 sector show the framework can accommodate high arrival demand without mid-air holding in TRACON airspace. The optimized behavior matches ATC practice, where speed control is used for small delay absorption and transitioning to geometric path stretching when larger spacing is needed. Across tested scenarios with realistic arrival rates, the objective structure strictly enforced safety, while maintaining runway throughput.

However, the current formulation stays with the FCFS rule. Under FCFS, trailing aircraft are forced into maximal path extensions to accommodate the sequence, even when a position swap can yield a more efficient global solution. Future work will extend this framework to include combinatorial sequence optimization, allowing for controlled position swaps within the rolling horizon. Additionally, we aim to incorporate stochastic wind uncertainty into the trajectory prediction model to further enhance the robustness of the generated schedules.





\clearpage
\bibliographystyle{scsproc}
\bibliography{ref}




\end{document}